\documentclass[a4paper,oneside]{amsart}
\usepackage[utf8]{inputenc}
\usepackage{lmodern,mathtools,amsmath,amsfonts,amssymb,url,xcolor,tikz,wrapfig}
\usepackage[colorlinks,urlcolor=blue,citecolor=green!50!black]{hyperref}
\usetikzlibrary{calc,angles,intersections,decorations.markings}

\makeatletter
\newcommand\wrapfigstop{\WF@finale}
\makeatother

\title{On angular measures in axiomatic Euclidean planar geometry}

\newcommand\mailto[1]{\href{mailto:#1}{#1}}

\author[M. Gr\"otschel]{Martin Gr\"otschel}
\address[M. Gr\"otschel]{Berlin-Brandenburg Academy of Sciences and Humanities\\
  J{\"a}gerstra{\ss}e 22/23\\
  10117 Berlin\\Germany}
\email{\mailto{groetschel@bbaw.de}}

\author[H. Hanche-Olsen]{Harald Hanche-Olsen}
\address[H. Hanche-Olsen]{Department of Mathematical Sciences\\
  NTNU Norwegian University of Science and Technology\\
  NO-7491 Trondheim\\ Norway}
\email{\mailto{harald.hanche-olsen@ntnu.no}}
\urladdr{\url{https://www.ntnu.edu/employees/hanche}}

\author[H. Holden]{Helge Holden}
\address[H. Holden]{Department of Mathematical Sciences\\
  NTNU Norwegian University of Science and Technology\\
  NO-7491 Trondheim\\ Norway}
\email{\mailto{helge.holden@ntnu.no}}
\urladdr{\url{https://www.ntnu.edu/employees/holden}}

\author[M. P. Krystek]{Michael P. Krystek}
\address[M. P. Krystek]{Physikalisch-Technische Bundesanstalt\\
  Bundesallee 100\\ D-38116 Braunschweig\\
  Germany}
\email{\mailto{Michael.Krystek@ptb.de}}

\newenvironment{axioms}
  {
    \begin{list}{X}{
        \setlength\leftmargin{2.5em}
        \setlength\labelsep{0.5em}
        \setlength\labelwidth{2em}
        }}
  {\end{list}}

\newenvironment{theorems}{\begin{axioms}\itshape}{\end{axioms}}

\newenvironment{singledashed}
  {\begin{scope}[decoration={
        markings,
        mark=at position 0.5 with {\draw[thin] (0pt,-2pt)--(0pt,2pt);}}]
      \newcommand\dpath{\path[postaction={decorate}]}
      \newcommand\ddraw{\draw[postaction={decorate}]}
    }
  {\end{scope}}

\newenvironment{doubledashed}
  {\begin{scope}[decoration={
        markings,
        mark=at position 0.5 with
        {\draw[thin] (-1pt,-2pt)--(-1pt,2pt) (1pt,-2pt)--(1pt,2pt);}}]
      \newcommand\dpath{\path[postaction={decorate}]}
      \newcommand\ddraw{\draw[postaction={decorate}]}
    }
  {\end{scope}}

\newenvironment{tripledashed}
  {\begin{scope}[decoration={
        markings,
        mark=at position 0.5 with
        {\draw[thin]
          (-2pt,-2pt)--(-2pt,2pt) (0pt,-2pt)--(0pt,2pt) (2pt,-2pt)--(2pt,2pt);}}]
      \newcommand\dpath{\path[postaction={decorate}]}
      \newcommand\ddraw{\draw[postaction={decorate}]}
    }
  {\end{scope}}

\newcommand{\ax}[1]{\textsf{#1}}

\DeclarePairedDelimiter\Set\{\}
\DeclarePairedDelimiterX\SetWhere[2]\{\}{\,#1\,\delimsize\vert\,\mathopen{}#2\,}

\newcommand\NN{\mathbb{N}}
\newcommand\QQ{\mathbb{Q}}
\newcommand\RR{\mathbb{R}}
\newcommand\Quant{\mathfrak{Q}}
\newcommand\Dimen{\mathfrak{D}}
\newcommand\Units{\mathfrak{U}}
\newcommand\Dim{\mathsf}
\newcommand\seg[1]{\overline{\mathit{#1}}}
\newcommand\len[1]{\langle\seg{#1}\rangle}
\newcommand\lin[1]{\mathit{#1}}
\newcommand\ray[1]{\overrightarrow{\mathit{#1}}}
\newcommand\ang[1]{\angle\mathit{#1}}
\newcommand\angm[1]{\langle\angle{\mathit{#1}}\rangle}
\newcommand\tri[1]{\mathit{#1}}
\newcommand\rad{\mathop{}\!\mathrm{rad}}

\newsavebox\gbox

\expandafter\providecommand\csname subjclassname@2020\endcsname
  {\textup{2020} Mathematics Subject Classification}


\date{\today}
\begin{document}

\begin{abstract}
  We address the issue of angular measure,
  which is a contested issue for the International System of Units (SI).
  We provide a mathematically rigorous and axiomatic presentation of angular measure
  that leads to the traditional way of measuring a plane angle
  subtended by a circular arc as
  the length of the arc divided by the radius of the arc, a scalar quantity.
  We distinguish between the \emph{angular magnitude},
  defined in terms of congruence classes of angles,
  and the (numerical) \emph{angular measure}
  that can be assigned to each congruence class
  in such a way that, e.g., the right angle has the numerical value $\frac\pi2$.
  We argue that angles are intrinsically different from lengths,
  as there are angles of special significance
  (such as the right angle, or the straight angle),
  while there is no distinguished length in Euclidean geometry.
  This is further underlined by the observation that,
  while units such as the metre and kilogram have been refined over time
  due to advances in metrology,
  no such refinement of the radian is conceivable.
  It is a mathematically defined unit, set in stone for eternity.
  We conclude that angular measures are numbers, and the current
  definition in SI should remain unaltered.
\end{abstract}

\maketitle

\bigskip
\begingroup\raggedleft\itshape\small
What’s in a name? That which we call a rose,\\
By any other word would smell as sweet.\\
\upshape – William Shakespeare
\par\endgroup

\section{Background}

Our motivation for the present study of a very classical question
is the ongoing discussion in the metrology community
regarding angular measure,
and, in particular,
whether one should associate a dimension to the angular measure.

The current version of SI states that angular measure is,
in the recent parlance of metrology, of dimension number.
(The term \emph{dimensionless} is still often used instead,
but is deprecated.)
A proposal has been raised to add the radian as an eighth base unit.
The purpose of the present paper and the companion paper \cite{comment}
is to show why we disagree with this proposal.
We discuss our approach in the context of
the ongoing discourse in the metrology community in \cite{comment}.
Here, we go deeper into the technical side of the argument.

To further explain this setting, we refer to \cite{Krystek}.
The set of measurable physical quantities $\Quant$
is an abelian group, with group operation written multiplicatively,
containing a copy (via an embedding $\iota$)
of the positive real numbers $\RR^{+}$.
(For the present discussion,
we ignore the possibility of negative quantities.)
The quotient group $\Dimen=\Quant/\RR^{+}$
is called the group of \emph{quantity dimensions}.
It is a free Abelian group
on finitely many generators, the \emph{base dimensions}.
Currently, the SI has seven base dimensions:
Length $\Dim{L}$,
mass $\Dim{M}$,
time $\Dim{T}$,
electric current $\Dim{I}$,
thermodynamic temperature $\Dim{\Theta}$,
amount of substance $\Dim{N}$, and
luminous intensity $\Dim{J}$.
The controversy concerns a suggested eighth base dimension
for angles.
The quotient map from $\Quant$ to $\Dimen$ is written $\delta$.
Thus we have the exact sequence
\begin{equation*}
\RR^{+} \stackrel{\iota}{\hookrightarrow} \Quant \stackrel{\delta}{\twoheadrightarrow} \Dimen.
\end{equation*}
Selecting a coherent set of units
yields a subgroup $\Units$ of $\Quant$
so that $\delta$ maps $\Units$ isomorphically onto $\Dimen$.
As a result, we have the isomorphism $\Quant \simeq \RR^{+}\times\Units$.
To illustrate these concepts, the acceleration of gravity $g$
has quantity dimension $\delta(g)=\Dim{M}\Dim{T}^{-2}$,
and $g \approx 9.8\,\mathrm{m}\,\mathrm{s}^{-2}$
corresponding to the pair
$(9.8,\mathrm{m}\,\mathrm{s}^{-2}) \in \RR^{+}\times\Units$.
Writing $\Dim{Z}$ for the identity element of $\Quant$,
we have $\iota(\RR^{+})=\delta^{-1}(\Dim{Z})$.
The quantities in this set correspond to pure numbers.
They are the quantities of dimension number.

Here we argue from a mathematical point of view
that angular measure cannot become a base quantity.
In support of this opinion,
we carefully revisit the concept of angle and angular measure
starting from classical Euclidean geometry
as recast in modern form by Hilbert.
We introduce (abstract) angular magnitude
as congruence classes of angles,
and associate an angular measure – being a real number –
to each angular magnitude.

As is evident from the present mathematical analysis,
the traditional way of measuring
a plane angle subtended by a circular arc,
is in the axiomatic approach to take the supremum
of sums of ratios of straight line length segments
of decreasing length.
The inevitable conclusion is that angular measures are pure numbers.

\section{Introduction}

Our goal here is to offer a detailed presentation of a mathematicians' view on
the question of angular measure,
with a focus on the mathematical concept of angles
rather than their physical manifestation.
In the following we have chosen to pursue the axiomatic approach
introduced in Euclid's \textit{Elements} \cite{euclid}
as made rigorous by D. Hilbert \cite{hilbert}.
Using the ingenious method of Archimedes \cite{heath},
we finally obtain the traditional way of measuring a plane angle
subtended by a circular arc
as the length of the arc divided by the radius of the arc, a scalar quantity.
For simplicity of presentation,
we limit our discussion to \emph{planar} geometry
throughout the paper.

\medskip
The axiomatization of planar geometry
as laid out in Euclid's \emph{Elements}
is a pillar in the development of mathematics.
In 1899, David Hilbert
gave a modern axiomatic formulation of Euclidean geometry based
on the  two millennia of mathematical progress since Euclid.
See \cite{gray} for a historical discussion.
We follow the lucid presentation due to Hartshorne \cite{hartshorne},
with a few twists of our own.

Our focus here is to give an axiomatic presentation
of angular measure in the Euclidean plane,
leading to the familiar definition of
the measure of a plane angle subtended by a circular arc
as the length of the arc divided by the radius of the arc,
a pure number.

Planar Euclidean geometry is well-known;
however, an axiomatic presentation requires a certain care
to develop the tools in a complete and consistent manner.
To ease the task for the reader
we give a rather detailed and complete presentation,
starting with the undefined notions of points and lines,
based on Hilbert's axioms as presented in \cite{hartshorne}.
These axioms include a notion of \emph{betweenness},
allowing us to define a line segment
as the set of points between two points on a line.
Furthermore, the axioms include a notion of congruence between line segments.
Thus, we can define the (abstract, or geometric) \emph{length}
of a line segment to be the congruence class of that line segments.
We can associate a real number to a pair of such lengths,
which we may think of as their length ratio.
We could assign some arbitrary length the rôle of unit length,
thus allowing to measure any line segment using real numbers.
However, we choose not to do so, staying with ratios instead.

Before presenting our approach precisely, we provide an informal preview.
An angle is defined as the union of two (distinct and non-opposite) rays
(half-lines with a direction)
originating from a common point (denoted apex).
Note that this excludes the zero angle and the straight angle.
A congruence relation between angles is introduced axiomatically.
We introduce the addition of (congruence classes of) angles
provided their sum is less than the straight angle
(more precisely, if each angle is less than the supplementary angle of the other).
The extension to angles of arbitrary magnitude
is essentially a book keeping issue.

Next comes the definition of the size of an angle.
While lengths scale indiscriminately,
angles are different in the sense that the right angle and the straight angle
are distinguished any way you measure them.
Here we follow the celebrated approach due to Archimedes
in his approximation of
the ratio of the circumference of a circle to its diameter
(the symbol $\pi$ for this ratio was introduced by William Jones in 1706,
and later popularised by Euler).
Archimedes' method consisted in
approximating the circle
by inscribed and circumscribed regular polygons of high order.
Using a regular 96-sided polygon, he arrived at the estimates
$3\frac{10}{71}<\pi<3\frac17$ \cite[pp.\ 93–98]{heath}.

We define the measure of an angle
as the supremum of the sum of length ratios
of the polygonal lines approximating
the circular arc subtending the angle
and the radius.
This notion of angle measure is additive
and gives the measure of $\frac\pi2$ for a right angle.
We note that the measure of an angle is a pure number,
arising as it does from sums of length ratios.

We distinguish between the concept of \emph{angular magnitude}
defined as a congruence class of angles,
and the  \emph{angular measure},
assigned to each congruence class.
More specifically, each congruence class consists of angles characterized
in such a way that they are ``of the same size''.
As we will argue, this is a function of length \emph{ratios}, and
thus independent of any length scale.
For each angular magnitude $\alpha$ we associate a numerical value  $\vartheta(\alpha)$,
which we can write as usual  $\vartheta(\alpha)=s/r$,
as the ratio between the arc length and the radius.
This will lead to the classical result
that a right angle has the numerical value $\frac\pi2$.
The radian is defined as
the angular magnitude for which $\vartheta(\rad)=1$.

The common conflation of identifying $\vartheta(\alpha)$ and $\alpha$
  appears to be the main source of much confusion regarding angular measure.
  In practical computations and measurement,
  this does not cause any problems,
  but they are conceptually different.
  In the present paper,
  we are only concerned with the theoretical aspects
  of angular measures,
  not with the considerable challenges
  associated with practical, accurate measurement of angles.
  Our goal is to participate from a mathematical point of view in
  the ongoing discussion regarding
  a base unit for angular measures.
  In light of our findings,
  we argue that it is neither desirable nor reasonable
  to add the radian as a base unit in the SI.

\section{Points, lines, and linear segments}
Like all axiomatic systems, geometry is based on a number of undefined terms.
In our case, we begin with a set whose elements are called \emph{points}
and another set whose elements are called \emph{lines}.
As we restrict our attention to \emph{planar} geometry,
we do not need the extra concept of a plane.
We will introduce, step by step,
a system of axioms that the points and lines have to satisfy.
Of particular importance are the notions of
\emph{incidence}, which is a relation between a point and a line,
\emph{betweenness}, which is a relation between three collinear points,
and \emph{congruence}  $\cong$, which is an equivalence relation\footnote
{An \emph{equivalence relation} $\cong$ between two objects $a$ and $b$
  is a relation that is
  reflexive ($a \cong a$),
  symmetric (if $a \cong b$, then $b \cong a$)
  and transitive (if $a \cong b$ and $b \cong c$, then $a \cong c$).
  The \emph{equivalence class} of an object $a$
  is the set of all objects that are equivalent to $a$.}
between line segments or between angles (to be defined later).

\savebox\gbox
{\begin{tikzpicture}[x=5mm,y=5mm]
  \def\segm(#1)#2:#3(#4){($(#1)!#2!(#4)$)--($(#1)!#3!(#4)$)}
  \begin{scope}[radius=0.5mm]
    \fill (0,2) coordinate (A) node[below]{\small$A$} circle;
    \fill (3,4) coordinate (B) node[right]{\small$B$} circle;
    \fill (4,1) coordinate (C) node[above right]{\small$C$} circle;
    \draw \segm(A)-0.2:1.2(B);
    \draw \segm(B)-0.2:1.2(C);
    \draw \segm(C)-0.2:1.2(A);
  \end{scope}
\end{tikzpicture}}
\begin{wrapfigure}[4]{r}{\wd\gbox}
  \vspace*{-\baselineskip}
  \usebox\gbox
\end{wrapfigure}

\medskip\noindent
The \emph{incidence axioms} are:
\begin{axioms}
\item[I1] Any two distinct points are incident with exactly one line.
\item[I2] Every line is incident with at least two distinct points.
\item[I3] There exist three noncollinear points.
\end{axioms}
It follows from \ax{I1} and \ax{I2} that any line is given
by the set of points incident with it.
Thus we can, and shall, identify a line
with its set of incident points.
Rather than the cumbersome “incident with” we use commonly understood
terms such as points lying on a line, a line passing through a point, etc.
We shall write $\lin{AB}$ for the unique line
through distinct points $A$, $B$.

\savebox\gbox
{\begin{tikzpicture}[x=5mm,y=5mm,radius=0.5mm]
    \def\segm(#1)#2:#3(#4){($(#1)!#2!(#4)$)--($(#1)!#3!(#4)$)}
    \fill (0,0) coordinate (A) node[below]{\small$A$} circle;
    \fill (2,0.6666) coordinate (B) node[below]{\small$B$} circle;
    \fill (3,1) coordinate (C) node[below]{\small$C$} circle;
    \draw \segm(A)-0.2:1.2(C);
    \fill (0,-2) coordinate (D) node[below]{\small$D$} circle;
    \fill (3,-1) coordinate (E) node[below]{\small$E$} circle;
    \draw[dotted] \segm(D)-0.2:1.2(E);
    \draw[very thick] (D)--(E) node[pos=0.5, below]{\small$\seg{DE}$};
  \end{tikzpicture}}
\begin{wrapfigure}{r}{\wd\gbox}
  \vspace*{-2\baselineskip}
  \usebox\gbox
\end{wrapfigure}
\medskip\noindent
For the \emph{betweenness axioms} \ax{B1}–\ax{B4},
we need some notation and a definition.
The betweenness relation “\!$B$ is between $A$ and $C$”
is written $A*B*C$.
The \emph{line segment} between two distinct points $D$ and $E$ is the set
consisting of $D$, $E$, and all points between them, and is denoted by
$\seg{DE}$.

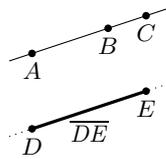
\wrapfigstop
\begin{axioms}
\item[B1]
  If $A*B*C$ then $A$, $B$, $C$ are distinct points on a line,
  and also $C*B*A$.
\item[B2] For any two distinct points $A$ and $B$,
  there exists a point $C$ such that $A*B*C$.
\item[B3] Given three distinct points on a line,
  exactly one of them is between the other two.
\item[B4] \emph{Pasch's axiom}:\footnote{%
    Strictly speaking, Pasch stated the axiom for three-dimensional geometry,
    with the added assumption that $l$ lies the plane containing $A$, $B$, and $C$.
    The present version, in contrast, restricts the dimensions to two.}
  If $A$, $B$, $C$ be three non-collinear points
  and a line $l$ contains none of them,
  but $l$ contains a point in $\seg{AB}$,
  then $l$ contains a point in $\seg{AC}\cup\seg{BC}$.
\end{axioms}

\savebox\gbox
{\begin{tikzpicture}[x=5mm,y=5mm]
  \begin{scope}[radius=0.5mm]
    \def\segm(#1)#2:#3(#4){($(#1)!#2!(#4)$)--($(#1)!#3!(#4)$)}
    \fill (0,2) coordinate (A) node[below]{\small$A$} circle;
    \fill (3,4) coordinate (B) node[right]{\small$B$} circle;
    \fill (4,1) coordinate (C) node[right]{\small$C$} circle;
    \draw[very thick] (A)--(B)--(C)--cycle;
    \fill ($(A)!0.6!(B)$) coordinate (AB) circle ;
    \coordinate (BC) at ($(B)!0.55!(C)$);
    \draw \segm(AB)-0.5:1.6(BC) node[pos=0,left]{\small$l$};
    \filldraw[fill=white] (BC) circle;
  \end{scope}
\end{tikzpicture}}
\begin{wrapfigure}{r}{\wd\gbox}
  \centering
  \usebox\gbox\par\footnotesize Pasch's axiom
\end{wrapfigure}

\noindent
While axiom \ax{I3} makes our geometry \emph{at least} two-dimensional,
Pasch's axiom \ax{B4} in effect makes it \emph{at most} two-dimensional.
One aspect of this is that
the set of points not on $l$
is divided into two nonempty disjoint subsets,
so that a line segment between two points in one subset
does not intersect $l$,
while a line segment between a point in one subset
and a point in the other does intersect $l$.
The two sets are called the two \emph{sides} (or half planes) of $l$,
and we are thus allowed to use phrases like
“\!$A$ and $B$ lie on the same side of $l$” or
“\!$A$ and $B$ lie on opposite sides of $l$”.

Axioms \ax{B1}–\ax{B3} imply that the set of points on any line $l$
can be given a total order $\prec$ so that for any distinct $A$, $B$, $C$
on $l$, $A*B*C$ if and only if either $A \prec B \prec C$ or $C \prec B \prec A$.
Moreover, the order is unique up to reversal.
For this reason, the axiomatic system
given by the incidence axioms \ax{I1}–\ax{I3}
and the betweenness axiom \ax{B1}–\ax{B4}
is known as \emph{ordered geometry}.

\savebox\gbox
{\begin{tikzpicture}[x=5mm,y=5mm,radius=0.5mm]
    \fill (0,0) node[above]{\small{$A$}} circle;
    \fill (3,0) node[above]{\small{$B$}} circle;
    \draw[dotted] (-1,0)--(0,0);
    \draw[very thick] (0,0)--(4,0);
    \draw[very thick,dashed] (4,0)--(5.5,0);
  \end{tikzpicture}}
\begin{wrapfigure}[2]{r}{\wd\gbox}
  \vspace*{-\baselineskip}
  \usebox\gbox
\end{wrapfigure}

The \emph{ray} $\ray{AB}$ consists of
all points $C \in \lin{AB}$ so that $C*A*B$ is \emph{not} true.
We say the ray \emph{originates} at $A$,
or that $A$ is the \emph{origin} of the ray.

\medskip\noindent
The \emph{congruence axioms for line segments} are:

\begin{axioms}
\item[C1]
  Given a line segment $\seg{AB}$ and a ray $r$ originating at $C$,
  there is a unique $D \in r$ so that $\seg{AB} \cong \seg{CD}$.
\item[C2]
  Congruence is an equivalence relation on line segments.
\item[C3]
  If $A*B*C$ and $D*E*F$
  and $\seg{AB} \cong \seg{DE}$
  and $\seg{BC} \cong \seg{EF}$,
  then $\seg{AC} \cong \seg{DF}$.
\end{axioms}

\savebox\gbox
{\begin{tikzpicture}[x=5mm,y=5mm,radius=0.5mm]
    \fill (0,0) coordinate(A) node[above]{\small{$A$}} circle;
    \fill (2,0) coordinate(B) node[above]{\small{$B$}} circle;
    \fill (-1,-1) coordinate(C) node[above]{\small$C$} circle;
    \draw[very thick] (A)--(B);
    \draw (C)-- +(-20:4) node[above left]{\small$r$};
    \draw[very thick] (C)-- +(-20:2) circle node[above]{\small$D$};
  \end{tikzpicture}}
\begin{wrapfigure}{r}{\wd\gbox}
  \centering
  \vspace*{-\baselineskip}
  \usebox\gbox\par\footnotesize Axiom \ax{C1}
\end{wrapfigure}

\noindent
As we like to think of line segments being congruent
if and only if they have the same length,
thanks to \ax{C2} we can turn this into a definition,
and define the \emph{length} $\len{AB}$
of a segment $\seg{AB}$ to be its congruence class.
Then axiom \ax{C1} states that the points on a ray
(except for the origin of the ray)
are in one-to-one correspondence with the set of lengths.
It follows from the axioms that no line segment
is congruent to a proper subsegment.
Thus we can define a total order on the set of lengths,
so that the length of a line segment is greater than that of any of its subsegments,
and the sum of two lengths is greater than either summand.

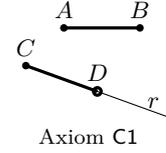
\savebox\gbox
{\begin{tikzpicture}[x=5mm,y=5mm,radius=0.5mm]
    \fill (0,0) coordinate(A) node[above]{\small{$A$}} circle;
    \fill (2,0) coordinate(B) node[above]{\small{$B$}} circle;
    \fill (3,0) coordinate(C) node[above]{\small{$C$}} circle;
    \fill (0,-1) coordinate(D) node[below]{\small{$D$}} circle;
    \fill (D) +(-20:2) coordinate(E) node[below]{\small{$E$}} circle;
    \fill (D) +(-20:3) coordinate(F) node[below]{\small{$F$}} circle;
    \draw[very thick] (A)--(C) (D)--(F);
  \end{tikzpicture}}
\begin{wrapfigure}{r}{\wd\gbox}
  \centering
  \vspace*{-\baselineskip}
  \usebox\gbox\par\footnotesize Axiom \ax{C3}
\end{wrapfigure}

We can \emph{add} lengths by saying
$\len{AB}+\len{BC}=\len{AC}$ whenever $A*B*C$;
axiom \ax{C3} states that this is indeed well defined.
Thus the set of lengths of line segments
becomes an abelian semigroup.\footnote
{This semigroup lacks a neutral element – a zero,
since we do not include degenerate “line segments”
consisting of a single point.
This is easily remedied if desired.}
In particular, we can define any positive integer multiple of a length
by repeated addition.
So far, we can only add lengths,
but we cannot define a ratio of lengths,
except for special cases where the ratio will be a rational number:
If $m$ and $n$ are natural numbers with $m\len{AB}=n\len{CD}$,
it makes sense to say that $\len{AB}/\len{CD}=n/m$.

This is a good place to introduce the axioms of Archimedes and Dedekind.
Here we should note that the former does in fact follow from the latter.
However, we include it due to its importance in the discussion below.
\begin{axioms}
\item[A] (Archimedes' axiom) Given two line segments, some integer
  multiple of the first segment is greater than the second.
\item[D] (Dedekind's axiom)
  If a line is a disjoint union of two nonempty sets $A$ and $B$,
  and no point in $B$ lies between two points of $A$
  and vice versa,
  there exists a point $P$ so that
  $\Set{P} \cup A$ and $\Set{P} \cup B$ are opposite rays
  originating at $P$.
\end{axioms}

\savebox\gbox
{\begin{tikzpicture}[x=5mm,y=5mm,radius=0.5mm]
    \fill (0,0) coordinate(A) node[below]{\small{$A$}} circle;
    \fill (2,0) coordinate(B) node[below]{\small{$B$}} circle;
    \fill (B) +(80:2) coordinate(C) node[right]{\small{$C$}} circle;
    \fill (B) +(80:3) coordinate(D) node[left]{\small{$D$}} circle;
    \fill ($(A)!1.3!(D)$) coordinate(E) node[right]{\small{$E$}} circle;
    \draw[name path=AB, very thick] (A)--(B);
    \draw (B)--($(B)!1.4!(D)$) (A)--($(A)!1.2!(E)$);
    \draw[name path=CE] ($(E)!-0.5!(C)$)--($(E)!2.4!(C)$);
    \filldraw[fill=white, name intersections={of=AB and CE}]
    (intersection-1) circle;
  \end{tikzpicture}}
\begin{wrapfigure}{r}{\wd\gbox}
  \centering
  \vspace*{-\baselineskip}
  \usebox\gbox\par
\end{wrapfigure}

\noindent
Archimedes' axiom is essential
for avoiding the existence of infinitesimal or infinite lengths,
whereas Dedekind's axiom guards against the existence of
“point-sized holes” in a line.
Related is the fact that there is no smallest length;
i.e., between any given distinct points $A$ and $B$
another point can be found.
(A proof is briefly indicated in the figure.
$C$ is any point not on the line $\lin{AB}$,
then $D$ and $E$ are picked by \ax{B2},
and Pasch's axiom (\ax{B4}) is used to show that
$\lin{EC}$ must intersect $\seg{AB}$.)

Given two lengths $x$ and $y$,
we can define their \emph{ratio} $x/y \in (0,\infty)$
by
\begin{equation*}
  x/y = \sup \SetWhere{a/b}{a,b\in\NN, ay \le bx}.
\end{equation*}
Archimedes' axiom guarantees that this ratio is positive and finite,
whereas Dede\-kind's axiom (together with the non-existence of any
smallest length) implies that every positive real number
is a ratio of lengths as defined above.
Relations like $(x+y)/z=x/z+y/z$ and $x/z<y/z \Leftrightarrow x<y$ are easy to show.

Having defined ratios, we can now define real multiples of lengths:
If $a$ is a positive real number and $u$ is a length,
we can define $au$ to be the unique length $x$ so that $x/u=a$.
We could also fix such a length $u$, call it the unit length,
and assign the real number $x/u$ to any length $x$.
This is customarily done in elementary geometry,
but it is by no means necessary.
One could also assign special names to several such lengths,
and use multiples of the resulting “units” to specify lengths,
as one does when referring to physical space.

\section{Angles}

We now turn to the study and measure of \emph{angles}.
In this section, we take ordered geometry (incidence and betweenness axioms)
as given, along with the congruence axioms \ax{C1}–\ax{C3}.
Axioms \ax{A} and \ax{D} are not needed for now.

\savebox\gbox
{\begin{tikzpicture}[x=5mm,y=5mm,radius=0.5mm]
    \coordinate (A) at (0,0);
    \coordinate (B) at (2,0);
    \coordinate (C) at (2,3);
    \draw pic[fill=black!30,angle radius=3mm] {angle=B--A--C};
    \fill (A) node[below]{\small{$A$}} circle;
    \fill (B) node[below]{\small{$B$}} circle;
    \fill (C) node[above left=-2pt]{\small{$C$}} circle;
    \draw[very thick] ($(A)!1.4!(B)$)--(A)--($(A)!1.4!(C)$);
  \end{tikzpicture}}
\begin{wrapfigure}{r}{\wd\gbox}
  \usebox\gbox
  \vspace*{-\baselineskip}
\end{wrapfigure}

An \emph{angle} is the union of two distinct and non-opposite rays,
called the \emph{legs} of the angle,
originating from a common point (the \emph{apex} of the angle).
Given two such rays $r$ and $s$,
we may write $\angle{rs}$ for the corresponding angle ($r \cup s$),
but the notation $\ang{BAC}$ is more commonly employed,
where the two rays are $\ray{AB}$ and $\ray{AC}$.
Note that $\ang{BAC}=\ang{CAB}$:
Angles (like lines and segments) do not have a specific orientation.
They do, however, have an inside and an outside:
A point not on either leg of $\ang{BAC}$ is said to be \emph{inside} $\ang{BAC}$
if it lies on the same side of $\lin{AB}$ as $C$,
and on the same side of $\lin{AC}$ as $B$.
Otherwise, it is \emph{outside} the angle.

The \emph{triangle} $\tri{ABC}$,
with $A$, $B$, $C$ non-collinear,
is the union $\seg{AB}\cup\seg{BC}\cup\seg{CA}$
together with the given order of the corners $A$, $B$, $C$.
Thus $\tri{ABC}$ and $\tri{BAC}$ are not considered to be the same triangle,
although as points sets they are identical.
Triangles $\tri{ABC}$ and $\tri{DEF}$ are called \emph{congruent}
if corresponding sides and angles are congruent, i.e.,
$\seg{AB}\cong\seg{DE}$,
$\seg{BC}\cong\seg{EF}$,
$\seg{CA}\cong\seg{FD}$,
$\ang{BAC}\cong\ang{EDF}$,
$\ang{CBA}\cong\ang{FED}$, and
$\ang{ACB}\cong\ang{DFE}$.

\medskip\noindent
The first two \emph{congruence axioms for angles}
are exact analogues of the congruence axioms
\ax{C1} and \ax{C2} for segments.

\begin{axioms}
\item[C4]
  Given an angle $\ang{BAC}$
  and a ray $\ray{DF}$,
  there exists a unique ray $\ray{DE}$
  on a given side of $\lin{DF}$
  such that $\ang{BAC}\cong\ang{EDF}$.
\item[C5]
  Congruence is an equivalence relation on angles.
\item[SAS] (“Side–angle–side.”)
  Triangles $\tri{ABC}$ and $\tri{DEF}$
  are congruent, provided
  $\seg{AB}\cong\seg{DE}$, $\ang{BAC}\cong\ang{EDF}$,
  and $\seg{AC}\cong\seg{DF}$.
\end{axioms}

\noindent
The analogue of \ax{C3} for angle addition is an easy consequence of
\ax{SAS}, so we do not need to add it as an axiom:
\begin{theorems}
  \item[C6] \textup{Theorem.}
  If $C$ is an interior point of $\angle{BAD}$,
  and $G$ is an interior point of $\angle{FEH}$
  with $\angle{BAC}\cong\angle{FEG}$ and $\angle{CAD}\cong\angle{GEH}$,
  then $\angle{BAD}\cong\angle{FEH}$.
\end{theorems}

\savebox\gbox
{\begin{tikzpicture}[x=5mm,y=5mm,radius=0.5mm]
    \coordinate (A) at (0,0);
    \coordinate (B) at (2,0);
    \coordinate (C) at (2,3);
    \draw pic[fill=black!30,angle radius=3mm] {angle=B--A--C};
    \fill (A) node[below]{\small{$A$}} circle;
    \fill (B) node[below]{\small{$B$}} circle;
    \fill (C) node[above left=-2pt]{\small{$C$}} circle;
    \draw[very thick] (A)--(B)--(C)--cycle;
    \begin{singledashed}
      \dpath (A)--(B);
    \end{singledashed}
    \begin{doubledashed}
      \dpath (A)--(C);
    \end{doubledashed}
    \begin{tripledashed}
      \dpath (B)--(C);
    \end{tripledashed}
  \end{tikzpicture}}
\begin{wrapfigure}[6]{r}{\wd\gbox}
  \vspace*{-\baselineskip}
  \usebox\gbox
\end{wrapfigure}

\noindent
Here we pause to emphasize that \ax{SAS} is the main bridge
connecting the notion of lengths to the notion of angles.
To be specific, consider an angle $\ang{BAC}$.
Then \ax{SAS} states, among other things,
that $\len{BC}$ is uniquely determined by $\len{AB}$, $\len{AC}$,
and the congruence class of  $\ang{BAC}$.
In the opposite direction we find that
$\len{AB}$, $\len{AC}$, and $\len{BC}$
uniquely determine the congruence class of  $\ang{BAC}$.
This is a consequence of the “Three Sides” congruence theorem:
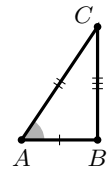
\begin{theorems}
\item[SSS] \textup{Theorem.}
  Given triangles $\tri{ABC}$ and $\tri{DEF}$,
  if $\seg{AB}\cong\seg{DE}$, $\seg{BC}\cong\seg{EF}$, and $\seg{CA}\cong\seg{FD}$,
  then the two triangles are congruent.
\end{theorems}
These results show that angular measures and length measures
are inextricably tied together,
a point we shall revisit later.

\medskip\noindent\emph{Remark.}
The incidence (\ax{I1}–\ax{I3}) and betweenness (\ax{B1}–\ax{B4}) axioms,
along with the congruence axioms \ax{C1}–\ax{C5} and \ax{SAS},
constitute an axiomatic system known as \emph{absolute geometry}.
Any model of absolute geometry, i.e., a set of lines and points
with the requisite relations satisfying the axioms,
is called a \emph{Hilbert plane}.
We call a Hilbert plane \emph{real}
if it satisfies the axioms of Archimedes (\ax{A}) and Dedekind (\ax{D})
in addition.
Here, “real” refers to the real numbers,
as these two axioms serve to ensure that length ratios
correspond to real numbers.

\section{The parallel postulate}
To complete our discussion of the axioms of geometry,
we introduce the \emph{parallel postulate}
according to Playfair
(Euclid's formulation was rather different).

Two lines
are said to be \emph{parallel}
if they have no point in common,
or else are the same line.

\begin{axioms}
  \item[P] (Playfair's axiom.)
  Given a line and a point,
  there is exactly one line through the given point
  parallel to the given line.
\end{axioms}

\noindent
\emph{Euclidean geometry} is absolute geometry
with Playfair's axiom added.
A \emph{Euclidean plane} is a model of Euclidean geometry,
i.e., a Hilbert plane satisfying axiom \ax{P}.
Finally, a \emph{real Euclidean plane}
is a Euclidean plane satisfying axioms \ax{A} and \ax{P}.

From now on, we shall assume all the axioms of a real Euclidean plane.
Non-Euclidean geometries have a number of peculiarities
that do not concern us here.

{\small
  A side remark: Similar triangles in a non-Euclidean geometry
  are congruent.
  In a Euclidean setting, on the other hand,
  even without axioms \ax{A} and \ax{D}
  we can define the ratio of lengths
  by an appeal to similar triangles.
  The resulting ratios are positive members of a
  Pythagorean ordered field,
  for example the real part of the algebraic closure of $\QQ$,
  or some non-standard model of $\RR$.
  \par}

\section{Angular magnitudes}

In the same way that \ax{C2} allows the definition of
lengths of intervals as congruence classes,
we use \ax{C5} to define the \emph{angular magnitude} of an angle
as its congruence class,
denoting by $\angm{BAC}$ the angular magnitude of $\ang{BAC}$.

\savebox\gbox
{\begin{tikzpicture}[x=5mm,y=5mm,radius=0.5mm]
    \coordinate (A) at (0,0);
    \coordinate (B) at (2,0);
    \coordinate (C) at (2,3);
    \coordinate (D) at (-2,2);
    \coordinate (D') at (-2,-1);
    \coordinate (E) at (-2,0);
    \fill (A) node[below]{\small{$A$}} circle;
    \fill (B) node[below]{\small{$B$}} circle;
    \fill (C) node[above left]{\small{$C$}} circle;
    \fill (D) node[below left]{\small{$D$}} circle;
    \fill (D') node[below right]{\small{$D'$}} circle;
    \fill (E) node[below]{\small{$E$}} circle;
    \draw[very thick]
    ($(A)!1.4!(B)$)--(A)--($(A)!1.4!(C)$)
    ($(A)!1.4!(C)$)--(A)--($(A)!1.4!(D)$)
    (A)--($(A)!1.4!(D')$);
    \draw (A)--($(A)!1.2!(E)$);
  \end{tikzpicture}}
\begin{wrapfigure}{r}{\wd\gbox}
  \vspace*{-\baselineskip}
  \centering
  \usebox\gbox
  \vspace*{-\baselineskip}
\end{wrapfigure}

Just as \ax{C3} with the help of \ax{C1}
lets us add arbitrary lengths,
we can use \ax{C6} with the help of \ax{C4} to define the sum
of angular magnitudes, with
the caveat that we cannot add angles
if the sum would be “too large”.
To make this more clear,
first note that we can order angles.
To do this, first move one of them
(by which we mean, replace it by a congruent angle)
so that the two angles to be compared have one leg in common,
with the other two legs on the same side of the common leg.
Unless the angles are congruent,
one will have its second leg inside the other angle,
and will be said to be smaller.
In the picture, $C$ is inside $\angle{BAD}$.
Hence $\angle{BAC}$ is smaller than $\angle{BAD}$,
and we write $\angm{BAC}<\angm{BAD}$.
The same picture serves also to define addition:
So long as $C$ is inside $\angle{BAD}$,
$\angm{BAC}+\angm{CAD}=\angm{BAD}$.
But clearly, two arbitrary angles cannot always
be arranged in this fashion, so their sum might not exist.

To be more precise, we can add angles if and only if
each angle is less than the supplementary angle of the other.
In the image, $\angm{BAC}+\angm{CAD}=\angm{BAD}$.
The supplementary angle of $\angle{BAC}$ is $\angle{CAE}$,
which is greater than $\angle{CAD}$.
It is, however, smaller than $\angle{CAD'}$,
and so the sum $\angm{BAC}+\angm{CAD'}$ is not defined.

To deal with the problem of adding “too large” angles,
the idea is to add (for the time being) a “fictional” angular magnitude $\varpi$
corresponding to the straight “angle” defined by opposite rays.\footnote
{The ad hoc notation $\varpi$ (a variant form of the letter $\pi$) is highly nonstandard.}
Technically, things become easier if we also adjoin a zero
to the set of angular magnitudes.
Then a generalized angular magnitude would be a formal sum $a\varpi+\varphi$,
where $a\ge0$ is an integer counting half turns,
and $\varphi$ is a proper angular magnitude (or zero).
The sum of $a\varpi+\varphi$ and $b\varpi+\psi$ would be one of
$(a+b)\varpi+(\varphi+\psi)$ or
$(a+b+1)\varpi+(\varphi+\psi-\varpi)$,
where we just have to give meaning to the term $\varphi+\psi-\varpi$
when $\varphi$ and $\psi$ cannot be properly added.
In the picture, $\angm{BAC}+\angm{CAD'}-\varpi=\angm{EAD'}$.
As a special case, the sum of a angular magnitude
and its supplement will be $\varpi$.
For future reference, note that a \emph{right} angle
is an angle congruent with its own supplement.
Its angular magnitude is $\varpi/2$.

It is a trivial, albeit quite tedious, book keeping exercise
to show that the set of generalized angular magnitudes
becomes an ordered additive semigroup
(in fact, a monoid, since we include the magnitude of the zero “angle”)
in analogy
with the semigroup of length measures.
We can define the ratio of angles
in analogy with how we defined ratios of lengths.

\savebox\gbox
{\begin{tikzpicture}[x=5mm,y=5mm,radius=0.5mm]
    \coordinate (A) at (0,0);
    \coordinate (B) at (3,0);
    \coordinate (C) at (50:3);

    \fill (A) node[left]{\small{$A$}} circle;
    \fill (B) node[below]{\small{$B$}} circle;
    \fill (C) node[above left]{\small{$C$}} circle;
    \draw pic[fill=black!30,angle radius=3mm] {angle=B--A--C};
    \draw (4,0)--(A)--(50:4);
    \draw[very thick] (B)--(A)--(C)--cycle;
    \begin{singledashed}
      \dpath (A)--(B);
      \dpath (A)--(C);
    \end{singledashed}
  \end{tikzpicture}}
\begin{wrapfigure}{r}{\wd\gbox}
  \vspace*{-\baselineskip}
  \centering
  \usebox\gbox
\end{wrapfigure}

We now proceed to the construction of angular measure.
We begin very naïvely,
relying on \ax{SAS} and \ax{SSS}
and, more generally,
the standard results on similar triangles,
which hold in Euclidean geometry – but not in non-Euclidean geometries –
and define
\begin{equation*}
  \sigma(\angm{BAC}) = \len{BC}/\len{AB}\qquad\text{provided $\seg{AB}\cong\seg{AC}$}.
\end{equation*}
In Euclidean geometry (i.e., satisfying \ax{P})
this ratio is well defined.
Moreover, it also makes intuitive sense:
An object of size $\len{BC}$ seen at a distance $\len{AB}$
subtends an angle measured by the ratio of the two lengths involved.

\savebox\gbox
{\begin{tikzpicture}[x=5mm,y=5mm,radius=0.5mm,line join=bevel]
    \coordinate (A) at (0,0);
    \coordinate (B) at (3,0);
    \coordinate (C) at (40:3);
    \coordinate (D) at (100:3);

    \fill (A) node[left]{\small{$A$}} circle;
    \fill (B) node[below]{\small{$B$}} circle;
    \fill (C) node[above left]{\small{$C$}} circle;
    \fill (D) node[above left]{\small{$D$}} circle;
    \draw (A)--(4,0) (A)--(40:4) (A)--(100:4);
    \draw[very thick] (B)--(C)--(D)--cycle;
    \begin{singledashed}
      \dpath (A)--(B);
      \dpath (A)--(C);
      \dpath (A)--(D);
    \end{singledashed}
  \end{tikzpicture}}
\begin{wrapfigure}{r}{\wd\gbox}
  \vspace*{-\baselineskip}
  \centering
  \usebox\gbox
\end{wrapfigure}

Indeed, $\sigma$ is an increasing function of the angular magnitude
as we have defined it.
However, it is not additive.
Rather, a simple application of the triangle inequality
reveals that it is \emph{subadditive}:
$\sigma(\alpha+\beta)<\sigma(\alpha)+\sigma(\beta)$,
as indicated in the picture,
with $\alpha = \angm{BAC}$ and $\beta = \angm{CAD}$.
We can remedy that by defining instead
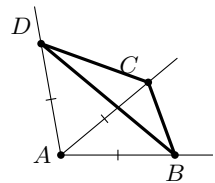
\begin{equation*}
  \vartheta(\alpha) =
  \sup \SetWhere[\bigg]{\sum_{i=1}^n \sigma(\beta_i)}{\sum_{i=1}^n \beta_i=\alpha}.
\end{equation*}
This is easily shown to be additive.
Briefly, first note that if
$\sum_i \beta_i = \alpha$ and $\sum_j \beta'_j = \alpha'$,
then $\sum_i \beta_i + \sum_j \beta'_j = \alpha + \alpha'$,
and so $\sum_i\sigma(\beta_i) + \sum_i\sigma(\beta'_i) \le \vartheta(\alpha + \alpha')$,
from which we get
$\vartheta(\alpha) + \vartheta(\alpha) \le \vartheta(\alpha + \alpha')$.
For the opposite inequality, if
$\sum_i \beta_i = \alpha + \alpha'$,
we may (if necessary) replace one of the $\beta_i$
by two angular magnitudes,
so that the angular magnitudes may be divided into two sets
summing to $\alpha$ and $\alpha'$, respectively.
Using the subadditivity of $\sigma$,
we find that this procedure increases
the value of $\sum_i\sigma(\beta_i)$,
and we get $\sum_i\sigma(\beta_i) \le \vartheta(\alpha) + \vartheta(\alpha')$,
so that $\vartheta(\alpha+\alpha') \le \vartheta(\alpha) + \vartheta(\alpha')$.

\savebox\gbox
{\begin{tikzpicture}[x=5mm,y=5mm,radius=0.5mm,line join=bevel]
    \coordinate (O) at (0,0);
    \coordinate (A) at (0:3);
    \coordinate (A') at (0:4);
    \coordinate (C) at (20:3);
    \coordinate (D) at (50:3);
    \coordinate (E) at (70:3);
    \coordinate (B) at (110:3);
    \coordinate (B') at (110:4);

    \draw (B')--(O)--(A');
    \draw[very thick] (A)--(C)--(D)--(E)--(B);

    \fill (O) node[left]{\small{$O$}} circle;
    \fill (A) node[below]{\small{$A$}} circle;
    \fill (B) node[left]{\small{$B$}} circle;

    \begin{singledashed}
      \dpath (O)--(A);
      \dpath (O)--(B);
      \ddraw (O)--(C);
      \ddraw (O)--(D);
      \ddraw (O)--(E);
    \end{singledashed}
    \draw[thick,radius=3] (A) arc[start angle=0, end angle=110];
  \end{tikzpicture}}
\begin{wrapfigure}{r}{\wd\gbox}
  \vspace*{-\baselineskip}
  \centering
  \usebox\gbox
\end{wrapfigure}

We can state our definition of $\vartheta(\alpha)$
in more geometric language as follows:
Given an angular magnitude $\alpha$ and angular magnitudes $\beta_i$
with $\sum_{i=1}^n \beta_i = \alpha$,
create an angle $\ang{AOB}$ with $\angm{AOB} = \alpha$,
pick a radius $r$,
and points $P_0=A$, $P_1$, \ldots, $P_n=B$
along the circular arc from $A$ to $B$
with $\angm{P_{i-1}OP_{i}}=\beta_i$ for $i=1$, \ldots, $n$.
Then
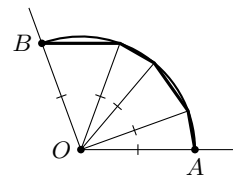
\begin{equation*}
  \sum_{i=1}^n \sigma(\beta_i)
  = \sum_{i=1}^n \frac{\len{P_iP_{i-1}}}{r}
  = \frac1r \sum_{i=1}^n \len{P_iP_{i-1}}
\end{equation*}
in which the sum on the right-hand side
is simply the length of the piecewise linear curve
passing from $A$ via the $P_i$ to $B$.
In the Cartesian plane $\RR^2$,
the length of a curve is defined to be
the supremum of the lengths of broken lines
formed by joining successive points along the curve.
We can employ the same definition in our more abstract setting,
concluding that $\vartheta(\alpha) = \ell/r$,
where $\ell$ is the length measure of the circular arc from $A$ to $B$:
\begin{equation*}
  \ell = \sup \sum_{i=1}^n \len{P_iP_{i-1}},
\end{equation*}
where the supremum is taken over all choices
of points $P_i$ picked successively along the arc.
The existence of the supremum in the set of length measures
is guaranteed by Dedekind's axiom (\ax{D}).
In this approach, the reader may recognize Archimedes' computation
of the circumference of a circle.
He used regular polygons,
approximating the circle both from the inside and the outside,
thus getting both a lower and an upper estimate.
But the idea is essentially the same.

\savebox\gbox
{\begin{tikzpicture}[x=5mm,y=5mm,radius=0.5mm,line join=bevel]
    \coordinate (O) at (0,0);
    \path[name path=HH] (-3,3)--(3,3);
    \coordinate (P) at (50:3);
    \path[name path=OP] (O)--(50:6);
    \path[name intersections={of=OP and HH,by=P'}];
    \coordinate (Q) at (75:3);
    \path[name path=OQ] (O)--(75:6);
    \path[name intersections={of=OQ and HH,by=Q'}];

    \pgfresetboundingbox

    \draw (P')--(O)--(Q');

    \fill (O) node[left]{\small{$O$}} circle;
    \fill (P) node[below]{\small{$P$}} circle;
    \fill (Q) node[below left]{\small{$Q$}} circle;
    \fill (P') node[above]{\small{$P'$}} circle;
    \fill (Q') node[above]{\small{$Q'$}} circle;

    \draw[thick,radius=3] (3,0) arc[start angle=0, end angle=180];
    \draw[thick] (3,0)--(3,3)--(-3,3)--(-3,0);
    \draw (P')--(O)--(Q');
    \draw[very thick] (P)--(Q) (P')--(Q');
    \begin{singledashed}
      \dpath (O)--(P);
      \dpath (O)--(Q);
    \end{singledashed}
  \end{tikzpicture}}
\begin{wrapfigure}{r}{\wd\gbox}
  \vspace*{-\baselineskip}
  \centering
  \usebox\gbox
\end{wrapfigure}

The one thing missing from the above discussion
is the fact that $\vartheta(\alpha)$
(equivalently, the length of the circular arc)
will in fact be finite.
The crucial observation here is that in the picture,
$\len{PQ} \le \len{P'Q'}$
(so long as $\lin{P'Q'}$ lies outside the circle).
Applying this to the segments of a broken line,
and using each of the three indicated sides
of the rectangle circumscribing the semicircle,
we quickly conclude that $\vartheta(\alpha)<4$ for any angular magnitude $\alpha$.

The ratio $\pi$ between the arc length of a semicircle and its radius
will be $\pi = \sup_\alpha \vartheta(\alpha)$,
the supremum taken over all proper angular magnitudes $\alpha$.
Thus we arrive at $\vartheta(\varpi)=\pi$.
If we perform the calculation in the standard Cartesian plane $\RR^2$,
we end up with the usual value for $\pi$,
\begin{equation*}
  \pi = \int_{-1}^{1} \frac{dx}{\sqrt{1-x^2}}.
\end{equation*}

\wrapfigstop

\section{Conclusion}

Let us repeat and look at the previous discussion from a general point of view.
Nobody will question that points and lines are geometric objects.
Angles, being the union of two rays with a common apex,
are geometric objects as well.
There is a practical need to associate numerical measures with geometric objects.
For lines, the natural way is to define a length unit (such as metre)
with which, for every line, the length of a line segment can be measured.
For angles, no matter how they are placed in the plane,
the natural way is to identify those angles that are “of the same size”
and define a measure that gives the same value to all angles of the same size.
We have captured this idea by the definition of congruence of angles and
have introduced the concept of angular magnitude in Euclidean plane geometry
as a congruence class of angles.
We have shown that the very notion of congruence of angles,
and hence the angular magnitudes, relies crucially on the concept of length.
However, since length units have no influence on angles,
we must conclude that angular measure must be considered
a function of length ratios.

Among the angular magnitudes we find $\varpi$, corresponding to the straight angle,
$\varpi/2$, corresponding to a right angle,
and the degree ${}^\circ=\varpi/180$.
To each angular magnitude $\alpha$ we have assigned an \emph{angular measure} $\vartheta(\alpha)$,
for which we can write in the conventional manner $\vartheta(\alpha)=s/r$.
In particular, $\vartheta(\varpi)=\pi$ and $\vartheta(1^\circ)=\pi/180$.
We can also define the radian as the angular magnitude for which $\vartheta(\rad)=1$.
We now have $\varpi=\pi\rad=180^\circ$;
in particular, we no longer need the temporary notation $\varpi$.

Note that the conventional notation $\alpha=s/r$ is, strictly speaking,
a category error, since
a angular magnitude is not a number.
It is, however, quite common to conflate the two concepts,
i.e., not to distinguish between $\alpha$ and $\vartheta(\alpha)$.
In the vast majority of cases this is harmless.

If we do conflate angular magnitudes with their numerical representation, however,
the equation $\vartheta(\rad)=1$ becomes $\rad=1$, which is the source of much confusion,
such as considering the radian
to be a derived unit which is equal to the number one.
Unfortunately, this statement also appears in the current SI brochure,
where moreover `$\rad$' is expressed by the quotient $\mathrm{m/m}$,
in order to emphasize that it \emph{is} a derived unit in the SI.
But these statements are not justified at all.

If any value associated with a magnitude is specified,
both the numerical value and the corresponding unit
must always be stated.
Angles are no exception.
In case of a semicircle, for example,
the value associated with the angular magnitude
shall be stated as $\pi\,\text{rad}$,
although $\pi$ would be sufficient from a mathematical point of view,
i.e., the “$\rad$” shall be added for clarification.
On the other hand, $c=\pi\,r$
must be written for the arc of a semicircle with radius $r$,
i.e., in this case it is necessary to omit the “$\rad$”,
because the angular measure has to be used here, which is a pure number.

We introduced the notion of
angular magnitude
and the conversion function $\vartheta$
only for the purpose of the present discussion.
However, requiring scientists and engineers
to maintain the distinction between angular magnitudes and their measure in radians
would impose an undue and totally unnecessary burden on them.
In particular, we do not propose the general use of our function $\vartheta$,
by  whatever name one would choose to give it.

At this point, we wish to make a point
regarding the fundamental nature of angles versus lengths
and other physical quantities.
Since the metre was introduced in 1793,
improvements in the science of metrology
has vastly increased the ability to measure lengths accurately,
in turn leading to the need to refine the very definition
of the metre in order to keep up with the technology.
No such claim can be made for angles.
In fact, even though we can certainly measure angles
much more accurately today than we could three centuries ago,
no conceivable technological advance can lead to a need
to refine the definition of the radian, or a right angle.
This simple observation supports the notion
that angle is a \emph{mathematical} concept
more than a topic of the physical sciences.
Mathematical objects do not require units for their measure,
as opposed to physical objects, which do.

Although the discussion here has been confined
to the angles of planar Euclidean geometry,
all conclusions apply equally to the concepts of
“angle of rotation” and “phase angle”,
which have not been discussed here
in order to concentrate on the essential points.

\end{document}